\documentclass[11pt,twoside]{amsart}
\usepackage{amsmath, amsthm, amscd, amsfonts, amssymb, graphicx, color}
\usepackage[bookmarksnumbered, colorlinks, plainpages]{hyperref}

\newtheorem{theorem}{Theorem}[section]
\newtheorem{lemma}[theorem]{Lemma}
\newtheorem{proposition}[theorem]{Proposition}
\newtheorem{corollary}[theorem]{Corollary}
\theoremstyle{definition}
\newtheorem{definition}[theorem]{Definition}

\theoremstyle{remark}

\numberwithin{equation}{section}

\begin{document}

\title[Some Results on Matricial Field $C^*$-Algebras]{Some Results on Matricial Field $C^*$-Algebras}
\author{Ali Ebadian}
\address{
Department of Mathematics \newline
\indent Urmia University, Urmia, Iran}
\email{ebadian.ali@gmail.com}
\author{Ali Jabbari}

\address{
Department of Mathematics \newline
\indent Urmia University, Urmia, Iran}
\email{jabbari\underline{ }al@yahoo.com}

\dedicatory{}



\begin{abstract}
In this paper, we consider Blackadar and Kirchberg’s MF algebras. We show that any  inner quasidiagonal C$^*$-algebra is MF algebra and we generalize Voiculescu's Representation Theorem for a special version of MF algebras. Moreover, we define a weak version of MF algebras namely matrical amenable (AM) algebras and prove some results related to this new notion. Finally, we consider real C$^*$-algebras and we show that a real C$^*$-algebra is MF if and only if its complexification is MF.
\end{abstract}
\maketitle
\section{Preliminaries and Introduction}
The notion of \emph{matricial field} (MF) C$^*$-algebras was introduced by Blackadar and Kirchberg in \cite{BK1996}. A separable  C$^*$-algebra is an MF algebra if it can be written as the inductive limit of a generalized inductive system of finite dimensional  C$^*$-algebras. This equivalent to that a separable  C$^*$-algebra $A$ is an MF algebra if and only if $A$ can be embedded as a C$^*$-subalgebra of $\left(\prod_n \mathbb{M}_{n_k}\right)/\left(\bigoplus_{n_k} \mathbb{M}_{n_k}\right)$, for a sequence of positive integers $n_k$, $k=1,2,\ldots$, where $\mathbb{M}_{n_k}$ is an $n\times n$ matrix over the complex space $\mathbb{C}$ \cite[Theorem 3.2.2]{BK1996}. Rainone and Schafhauser in \cite{RS2019} characterized MF algebras and gave interesting results related to the notion of MF algebras. Following Rainone and Schafhauser, a separable C$^*$-algebra $A$ is  called  (MF) if for every finite set $F\subseteq A$ and $\varepsilon>0$, there is an $n\geq1$ and $*$-linear map $\varphi:A\longrightarrow \mathbb{M}_n$ such that
\begin{itemize}
  \item[(i)] $\|\varphi(aa')-\varphi(a)\varphi(a')\|<\varepsilon$, and
  \item[(ii)] $\|\varphi(a)\|>\|a\|-\varepsilon$,
\end{itemize}
for all $a,a'\in F$.  A trace $\tau$ on $A$ is  called \emph{matricial field} (MF) if for every finite set $F\subseteq A$ and $\varepsilon>0$, there is an $n\geq1$ and $*$-linear map $\varphi:A\longrightarrow \mathbb{M}_n$ such that
\begin{itemize}
  \item[(i)] $\|\varphi(aa')-\varphi(a)\varphi(a')\|<\varepsilon$, and
  \item[(ii)] $|\tau(a)-\mathrm{tr}_n(\varphi(a))|<\varepsilon$,
\end{itemize}
for all $a,a'\in F$. MF algebras are a generalization of quasidiagonal (QD) C$^*$-algebras and indeed, for nuclear C$^*$-algebras this two notions, i.e., MF algebras and QD C$^*$-algebras are equivalent. A  unital $C^*$-algebra $A$ is called quasidiagonal, if, for every finite subset ${F}$ of $A$ and $\varepsilon>0$, there exist a matrix algebra $\mathbb{M}_k$ and a unital c.p. linear map $\varphi:A\longrightarrow \mathbb{M}_k$ such that
\begin{itemize}
  \item[(i)] $\|\varphi(ab)-\varphi(a)\varphi(b)\|<\varepsilon$, for all $a,b\in F$ and
  \item[(ii)] $|\ \|\varphi(a)\|-\|a\||<\varepsilon$, for all $a\in A$.
\end{itemize}

More connections between MF algebras and the classification of C$^*$-algebras can also be found in \cite{BK1996, BK2001, BK2011, HS2010}.

A  C$^*$-algebra $A$ is inner quasidiagonal if, for every $x_1, . . . , x_n \in A$
and $\varepsilon> 0$, there is a representation $\pi$ of $A$ on a Hilbert space $\mathcal{H}$ and a finite-rank
projection $P\in \pi(A)''$ such that $\|P\pi(x_j)-\pi(x_j)P\|<\varepsilon$ and $\|P\pi(x_j)P\|>\|x_j\|-\varepsilon$ for $1\leq j \leq n$ \cite{BK2001}.

Let $A$, $B$  and $D$ be  C$^*$-algebras with unital embeddings $\psi_A:D\longrightarrow A$ and $\psi_B:D\longrightarrow B$ the corresponding
full amalgamated free product C$^*$-algebra is the C$^*$-algebra $C$, equipped with unital embeddings $\sigma_A : A\longrightarrow C$ and $\sigma_B : B\longrightarrow C$ such that $\sigma_A\circ \psi_A = \sigma_B \circ \psi_B$, such that $C$ is generated by $\sigma_A(A)\cup \sigma_B(B)$ and satisfying the universal property that whenever $E$ is a C$^*$-algebra and $\pi_A : A\longrightarrow E$ and $\pi_B : B\longrightarrow E$ are $*$-homomorphisms satisfying $\pi_A\circ \psi_A = \pi_B \circ \psi_B$, there is a $*$-homomorphism $\pi : C \longrightarrow E$ such that $\pi\circ \sigma_A = \pi_A$ and $\pi\circ \sigma_B = \pi_B$. The full amalgamated free product C$^*$-algebra $C$ is commonly denoted by $A*_D B$. The full amalgamated free product of MF algebras is investigated in \cite{LS2012}.

Let $\mathcal{H}$ be a Hilbert space, then $B(\mathcal{H})$ denotes the set of bounded linear on $\mathcal{H}$. A family of elements $\{x_k\}_{k=1}^\infty$ in $B(\mathcal{H})$ converges to $x\in B(\mathcal{H})$ in *-SOT (*-strong operator topology) if and only if $x_k\to x$ in SOT $x_k^*\to x$ in SOT.   Let $\{A_k\}_{k=1}^\infty$ be a sequence of C$^*$-algebras such that $A_k\subseteq B(\mathcal{H}_k)$, where $\mathcal{H}_k$, for $k=1,2,\ldots$ is a Hilbert space. The direct product of $A_k$'s is defined as follows
$$\prod_kA_k=\left\{\langle a_k\rangle_{k=1}|\sup_k\|a_k\|<\infty, \quad a_k\in A_k\right\},$$
where the norm of any element $\langle a_k\rangle_{k=1}$ in $\prod_kA_k$ is defined by $\|\langle a_k\rangle_{k=1}\|=\sup_k\|a_k\|$. Moreover, the direct sum $A_k$'s is defined as follows
$$\bigoplus_kA_k=\left\{\langle a_k\rangle_{k=1}\in \prod_kA_k|\lim\sup_{k\to\infty}\|a_k\|=0, \quad a_k\in A_k\right\}.$$

The direct sum $\bigoplus_kA_k$ is a closed two sided of the C$^*$-algebra, so $\left(\prod_kA_k\right)/\left(\bigoplus_kA_k\right)$ becomes a C$^*$-algebra. A finite set $M\subseteq B(\mathcal{H})$ is quasidiagonal if there is an increasing sequence of finite-rank projections $\{P_n\}_{n=1}^\infty$ on $\mathcal{H}$ tending strongly to the identity such that $\|P_nT-TP_n\|\to0$ as $n\to\infty$. Following \cite{HS2010}, let $\mathbb{C}\langle X_1,\ldots,X_n\rangle$ be the set of all noncommutative polynomials in the indeterminates $X_1,\ldots,X_n, X_1^*,\ldots,X_n^*$. We need the following Lemmas in the next section, so we recall them as follow:
\begin{lemma}\label{lemi}\cite[Lemma 2.3]{HS2010}
  Suppose that $A$ is a unital C$^*$-algebra generated by a family of elements $x_1,\ldots , x_n$ in $A$. Then $A$ is MF if and only if for every faithful $*$-representation of $A$ on an infinite dimensional separable complex Hilbert space $\mathcal{H}$, $\pi:A\longrightarrow B(\mathcal{H})$ there is a family of elements $\left\{a_1^{(k)},\ldots, a_n^{(k)}\right\}\subseteq B(\mathcal{H})$ such that
\begin{itemize}
  \item[(i)] for each $k\geq 1$, $\left\{a_1^{(k)},\ldots, a_n^{(k)}\right\}$ is quasidiagonal.
  \item[(ii)] for any $P\in \mathbb{C}\langle X_1,\ldots,X_n\rangle$, $\lim_{k\to\infty}\left\|P\left(a_1^{(k)},\ldots, a_n^{(k)}\right)\right\|=\|P\left(x_1,\ldots,x_n\right)\|$.
  \item[(iii)] for every $1\leq i\leq n$, $a_i^{(k)}\to\pi(x_i)$ in *-SOT as $k\to\infty$.
\end{itemize}
\end{lemma}

\begin{lemma}\label{lemi1}\cite[Lemma 7.2.2]{BO2008}
Let $A$ be a unital C$^*$-algebra and let $pmq\in A$ be projections.
\begin{itemize}
  \item[(i)] If $\|p-q\|<1$, then there is a unitary element $u\in A$ with $uqu^*=p$ and $\|1_A-u\|<4\|p-q\|$.
  \item[(ii)] If $\|q-pq\|<1/4$, the there is a unitary element $u\in A$ such that $uqu^*\leq p$ and $\|1_A-u\|\leq 10\|q-pq\|$.
\end{itemize}
\end{lemma}

Let $\mathcal{Q}$ be the universal UHF algebra, $\omega$ be a free ultrafilter on $\mathbb{N}$ and let $\mathcal{Q}_\omega$ be the ultrapower of $\mathcal{Q}$ defined by
\begin{equation*}\label{up}
    \mathcal{Q}_\omega:=\ell^\infty(\mathcal{Q})/\{(a_n)\in\ell^\infty(\mathcal{Q})\ :\ \lim_{n\to\infty}\|a_n\|=0\}.
\end{equation*}

The unique trace $\tau_{\mathcal{Q}}$ on $\mathcal{Q}$ induces a trace $\tau_{\mathcal{Q}_\omega}$ on $\mathcal{Q}_\omega$ as follows:
$$\tau_{\mathcal{Q}_\omega}(x)=\lim_\omega\tau_{\mathcal{Q}}(x_n),\quad x\in \mathcal{Q}_\omega\ \text{is represented by}\ (x_n)\in\ell^\infty(\mathcal{Q}).$$

In the next section, we give some results on MF algebras and we generalize the Voiculescu's Representation Theorem that says that a C$^*$-algebra $A$ is QD if and only if $A$ has a faithful QD representation if and only if  every faithful unital essential representation of $A$ is QD. In section 3, we define MA algebras and MA races as generalizations of MF algebras and MF traces. We show that in case unital C$^*$-algebras MF traces and MA traces are coincide. Section 4 is related to real C$^*$-algebras and we investigate these algebras as MF algebras.
\section{Some results on MF algebras}
In this section, we give some results related to MF algebras. We commence with the following result that shows there exits a relationship between inner quasidiagonal C$^*$-algebras and  MF algebras.
\begin{proposition}\label{p1-1}
  Every separable inner quasidiagonal C$^*$-algebra is  MF.
\end{proposition}
\begin{proof}
  We do the same fashion as the proof of \cite[Lemma 4.1]{BK2011}. Let $A$ be inner quasidiagonal, then for finite subset $F$ of $A$
and $\varepsilon> 0$, there is a representation $\pi$ of $A$ on a Hilbert space $\mathcal{H}$ and a finite-rank
projection $P\in \pi(A)''$ such that $\|P\pi(a)-\pi(a)P\|<\varepsilon$ and $\|P\pi(a)P\|>\|a\|-\varepsilon$ for each $a\in A$. Let $Q_1, \ldots, Q_m$ and $R_1,\ldots, R_m$ be as in the proof of \cite[Lemma 4.1]{BK2011}. Then, we define $\rho=\pi|_{R_j\mathcal{H}}$ and $\varphi(a)=P_j\rho(a)P_j\in B(P_j\mathcal{H})\cong \mathbb{M}_n$, where $n = \dim(P_j\mathcal{H})$; for all $a\in A$. Thus, $\varphi:A\longrightarrow\mathbb{M}_n$ is a $*$-linear map. Then, according to definition of $\varphi$, $A$ is an MF algebra.
\end{proof}
\begin{corollary}
 A full amalgamated free product of two inner quasidiagonal separable C$^*$-algebras may not be inner quasidiagonal.
\end{corollary}
\begin{proof}
  Let $D=\mathbb{C}\oplus\mathbb{C}$ and consider the full amalgamated free product $\mathbb{M}_2(\mathbb{C})*_D\mathbb{M}_3(\mathbb{C})$ that is not MF, see \cite{LS2012}, for more details. Assume towards a contradiction $\mathbb{M}_2(\mathbb{C})*_D\mathbb{M}_3(\mathbb{C})$ is  inner quasidiagonal. Then by Proposition \ref{p1-1}, $\mathbb{M}_2(\mathbb{C})*_D\mathbb{M}_3(\mathbb{C})$  is MF, a contradiction.
\end{proof}

\begin{definition}
 Let $\mathcal{H}$ be an infinite dimensional Hilbert space and $M\subset B(\mathcal{H})$. We say that $M$ is a MF set if for each finite sets $F\subset M$, $H\subset \mathcal{H}$ and $\varepsilon>0$, there is a projection $P\in B(\mathcal{H})$ such that $\|PT-TP\|<\varepsilon$ and $\|Px\|>\|x\|-\varepsilon$, for all $T\in F$ and $x\in H$.
\end{definition}
\begin{lemma}\label{ll1}
  Let $\mathcal{H}$ be an infinite dimensional separable Hilbert space and $M\subset B(\mathcal{H})$ be a norm separable MF set. Then there is an increasing sequence of projections $P_1\leq P_2\leq\cdots$ such that $\|P_nT-TP_n\|\to0$  and $\|P_nx\|>\|x\|-\varepsilon$, for all $T\in M$ and $x\in \mathcal{H}$.
\end{lemma}
\begin{proof}
Let  $F\subset M$, $H\subset \mathcal{H}$ be a finite subsets, $\varepsilon>0$. Since $M$ is a MF set, there is a projection $\overline{P}\in B(\mathcal{H})$ such that $\|\overline{P}T-T\overline{P}\|<\varepsilon$ and $\|\overline{P}x\|>\|x\|-\varepsilon$, for all $T\in F$ and $x\in H$.
Let $Q$ be the orthogonal projection onto the span of $H$. Then one can find a larger finite set $\overline{H}$ of $Q\mathcal{H}$  such that $\|\overline{P}Q-Q\|<\varepsilon_1$, where $\varepsilon_1>0$ depends on $\varepsilon$. Then by  Lemma \ref{lemi1}, there is a unitary element $U\in B(\mathcal{H})$ and $0<c<\infty$ such that $Q\leq U\overline{P}U^*$ and $\|U-I\|<c\varepsilon_1$. Then, we have $\|\overline{P}-U\overline{P}U^*\|<2c\varepsilon_1$.

Now, we set $P_1=U\overline{P}U^*$. Then we have $\|P_1T-TP_1\|<\varepsilon$ and $\|P_1x\|>\|x\|-\varepsilon$, for all $T\in M$ and $x\in \mathcal{H}$. By continue this process we can build a sequence of projections $\{P_n\}_{n=1}^\infty$ such that satisfies $\|P_nT-TP_n\|\to0$, for all $T\in M$.
\end{proof}
\begin{definition}
Let $A$ be a C$^*$-algebra and $\mathcal{H}$ be an infinite dimensional separable Hilbert space. We call a representation $\pi:A\longrightarrow B(\mathcal{H})$ is a MF representation if $\pi(A)$ is a MF set of operators.
\end{definition}
By the following result we write Voiculescu's Representation Theorem for MF algebras.
\begin{theorem}
 Let $A$ be a unital separable C$^*$-algebra generated by $x_1,\ldots,x_n\in A$. Then the following statements are equivalent:
\begin{itemize}
  \item[(i)] $A$ is MF.
  \item[(ii)] $A$ has a faithful MF representation.
  \item[(iii)] every faithful unital essential representation of $A$ is MF.
\end{itemize}
\end{theorem}
\begin{proof}
  (i)$\to$(iii) Let $\pi:A\longrightarrow B(\mathcal{H})$ be a faithful unital essential representation, where $\mathcal{H}$ is a (infinite dimensional) separable Hilbert space. By Lemma \ref{lemi}, there is a family of elements $\left\{a_1^{(k)},\ldots,a_n^{(k)}\right\}_{k=1}^\infty\subseteq B(\mathcal{H})$ such that it is quasidiagonal and $a_i^{(k)}\to \pi(x_i)$ in *-SOT as $k\to\infty$. Thus, there is an increasing sequence of finite-rank projections $\{P_m\}_{m=1}^\infty$ on $\mathcal{H}$ tending strongly to the identity such that $\|P_ma_i^{(k)}-a_i^{(k)}P_m\|\to0$ as $m\to \infty$, for any $1\leq i\leq n$. This implies that
\begin{align*}
 \lim_{m\to\infty} \|P_m\pi(x_i)-\pi(x_i)P_m\| & =\lim_{m\to\infty}\left(\mathrm{*-SOT-}\lim_{k\to\infty}\|P_ma_i^{(k)}-a_i^{(k)}P_m\|\right) \\
   &= \mathrm{*-SOT-}\lim_{k\to\infty}\lim_{m\to\infty}\|P_ma_i^{(k)}-a_i^{(k)}P_m\|\to0.
\end{align*}

Hence, $\pi$ is MF.


For completing the proof it suffices we show that (ii) implies (i), because (iii)$\to$(ii) obviously holds. Let $\pi:A\longrightarrow B(\mathcal{H})$ be a faithful MF representation. Then by Lemma \ref{ll1}, there is an increasing sequence of projections $P_1\leq P_2\leq\cdots$ such that $\|P_nT-TP_n\|\to0$  and $\|P_nx\|>\|x\|-\varepsilon$, for all $T\in \pi(A)$ and $x\in \mathcal{H}$. Now, define $\varphi_n:A\longrightarrow P_nB(\mathcal{H})P_n\cong {M}_{k(n)}(\mathbb{C})$ by $\varphi_n(a)=P_n\pi(a)P_n$, for all $a\in A$. Then by easy calculations one can verify that $\varphi_n$'s satisfy in the conditions (i) and (ii) of MF algebras.
\end{proof}

\section{MA algebras}
In this section, we generalize MF algebras that we call them MA algebras and give some results related to these algebras.
\begin{definition}
 A trace $\tau$ on a separable C$^*$-algebra $A$ is  called \emph{matricial amenable} (MA) if for every finite set $F\subseteq A$ and $\varepsilon>0$, there is an $n\geq1$ and $*$-linear map $\varphi:A\longrightarrow \mathbb{M}_n$ such that
\begin{itemize}
  \item[(i)] $\|\varphi(aa')-\varphi(a)\varphi(a')\|_{2,\mathrm{tr}}<\varepsilon$, and
  \item[(ii)] $|\tau(a)-\mathrm{tr}_n(\varphi(a))|<\varepsilon$,
\end{itemize}
for all $a,a'\in F$, where $\mathrm{tr}$ is the unique normalised trace on the matrix algebra $\mathbb{M}_n$ and $\|\cdot\|_{2,\mathrm{tr}}$ is the induced norm on $\mathbb{M}_n$ that is $\|\alpha\|_{2,\mathrm{tr}}=\mathrm{tr}\left(\alpha^*\alpha\right)^\frac{1}{2}$.
\end{definition}

Note that if we can choose the above $\varphi$ from c.p.c maps, then the trace $\tau$ is called amenable.
\begin{lemma}\label{lem2.1}
  Let $A$ be a nonunital separable C$^*$-algebra and $\tau$ be a trace on $A$. Then there is a unique trace $\tilde{\tau}$ on $A^\#$ the unitization of $A$ such that $\tau$ is MF (MA) if and only if $\tilde{\tau}$ is MF (MA).
\end{lemma}
\begin{proof}
  By \cite[Proposition 3.5.10]{B2006}, there is a unique trace $\tilde{\tau}$ on $A^\#$ such that it is the extension of $\tau$. Let $\varphi:A\longrightarrow\mathbb{M}_n$ be a $*$-linear map, then we have the unital $*$-linear map $\tilde{\varphi}:A^\#\longrightarrow\mathbb{M}_n$ defined by $\tilde{\varphi}((a,\lambda))=\varphi(a)+\lambda 1_n$, where $1_n$ is the unit of $\mathbb{M}_n$. Then by easy calculations the proof will be complete.
\end{proof}
The proof of the following result is straightforward.
\begin{lemma}\label{lem2}
  Let $A$ be a separable C$^*$-algebra, $I$ be an ideal of $A$, $\tau$ is a trace on $A$ and $\tau'=\frac{1}{\|\tau|_I\|}\tau$ is a trace on $I$. If $\tau$  is  MF (MA), then $\tau'$ is also MF (MA).
\end{lemma}
\begin{definition}
  Let $A$ and $B_n$ be (unital) C$^*$-algebras where $n\in\mathbb{N}$ and $A$ is separable. A (unital) approximate morphism from $A$ into $(B_n)_{n\geq1}$ is a sequence of (unital) $*$-linear maps $(\varphi_n:A\longrightarrow B_n)_{n\geq1}$ satisfying
  $$\lim_{n\to\infty}\|\varphi_n(aa')-\varphi(a)\varphi(a')\|=0,$$
  for all $a,a'\in A$. Let $\tau\in \mathrm{T}(A)$ and $\tau_n\in\mathrm{T}(B_n)$, then the approximate morphism $(\varphi_n:A\longrightarrow B_n)_{n\geq1}$ is called trace-preserving if
  $$\lim_{n\to\infty}\tau_n(\varphi_n(a))=\tau(a),$$
  for all $a\in A$.
\end{definition}
\begin{lemma}\label{lem3}
Let $A$ be a separable \emph{(}unital\emph{)} C$^*$-algebra and $\tau$ be a trace on $A$. The following statements are equivalent:
\begin{itemize}
  \item[(i)] $\tau$ is MA.
  \item[(ii)] There is a trace-preserving approximate morphisms $(\varphi_n:A\longrightarrow\mathbb{M}_n)_{n\geq1}$.
  \item[(iii)] There is a trace-preserving  homomorphism $\phi:A\longrightarrow\prod_\omega\mathbb{M}_n$.
\end{itemize}
\end{lemma}
\begin{proof}
(i)$\to$(ii)  Similar to the proof of \cite[Proposition 2.2]{RS2019}, let $A_0$ be a countable, dense $*$-algebra over $\mathbb{Q}[i]$. Consider an increasing sequence $(F_n)_{n\in\mathbb{N}}\subseteq A$ of finite, self-adjoint subsets such that $A_0\subseteq \cup_{n\in\mathbb{N}}F_n$ and for any $a\in F_n$, there exists $b\in F_{n+1}$ and $\lambda\in\mathbb{C}$,
 \begin{equation}\label{eq1}
a^*a+\left(b+\lambda1_{A^\#}\right)^*\left(b+\lambda1_{A^\#}\right)=\|a\|\cdot 1_{A^\#}.
 \end{equation}

  Since $\tau$ is MA, for every $\varepsilon>0$ there is an $n\geq 1$ and a $*$-linear map $\phi_n:A\longrightarrow \mathbb{M}_n$ such that
  $$\|\phi_n(aa')-\phi_n(a)\phi_n(a')\|_{2,\mathrm{tr}}<\varepsilon_n\quad\text{and}\quad|\tau(a)-\mathrm{tr}_n(\phi_n(a))|<\varepsilon_n,$$
  for all $a,a'\in F_n$, where $\varepsilon_n>0$ is depend on $\varepsilon$. Fix $a\in A_0$, then there is a $k\geq1$ and $b\in F_{k+1}$ such that satisfying \eqref{eq1}. Then
  \begin{equation}\label{eq2}
   \phi_n(a^*a)+\phi_n\left(b^*b\right)+\lambda\phi_n\left(b^*\right)+\bar{\lambda}\phi_n\left(b\right)+|\lambda|^2\mathbb{I}_n=\|a\|\cdot \mathbb{I}_n.
  \end{equation}

  Then by \eqref{eq2}, we have
  \begin{align}\label{eq3}
  \nonumber
    \|\phi_n(a^*)\phi_n(a)-\left(\phi_n(b)+\lambda \mathbb{I}_n\right)^*\left(\phi_n(b)+\lambda\mathbb{I}_n\right)-\|a\|\cdot \mathbb{I}_n\|_{2,\mathrm{tr}}&\leq  \|\phi_n(a^*)\phi_n(a)- \phi_n(a^*a)\|_{2,\mathrm{tr}}\\
    \nonumber
    &\ + \|\phi_n(b^*)\phi_n(b)- \phi_n(b^*b)\|_{2,\mathrm{tr}}\\
    &<2\varepsilon_n,
  \end{align}
  where $n\geq k+1$. Thus, for such $n$, $\phi_n(a^*)\phi_n(a)\leq\left(\|a\|+2\varepsilon_n\right)\cdot\mathbb{I}_n$. Then
  $$\|\phi_n(a)\|_{2,\mathrm{tr}}=\mathrm{tr}(\phi_n(a^*)\phi_n(a))^\frac{1}{2}\leq \left(\|a\|+2\varepsilon_n\right)\mathrm{tr}\left(\mathbb{I}_n\right)^\frac{1}{2}=\|a\|+2\varepsilon_n,$$
  for $n\geq k+1$. Thus, sequence $(\phi_n(a))_n$, for $a\in A_0$ is bounded and moreover $\lim\sup_{n\to\infty}\|\phi_n(a)\|_{2,\mathrm{tr}}\leq\|a\|$, for $a\in A_0$.

   Let $\bigoplus_n\mathbb{M}_n$ be the norm closed ideal consisting of all sequences $(m_n)_n\subseteq\prod_n \mathbb{M}_n$ such that $\lim_n\|b_n\|_{2,\mathrm{tr}}=0$. Then similar to the proof of \cite[Proposition 2.2]{RS2019}, the maps $\phi_n$ induce a contractive $*$-homomorphism $\phi_0:A_0\longrightarrow\prod_n \mathbb{M}_n/\bigoplus_n\mathbb{M}_n$ which extends by continuity to a homomorphism $\phi:A\longrightarrow \prod_n \mathbb{M}_n/\bigoplus_n\mathbb{M}_n$ such that $\phi_0$ is isometric, $\phi$ is faithful and
    $\tau(a)=\lim_{n\to\infty}\mathrm{tr}_n(\phi(a))$.

    Suppose that $(\phi_n:A\longrightarrow \mathbb{M}_n)_{n\geq1}$ is a sequence of $*$-linear maps lifting $\phi$. Then
   $$\lim_{n\to\infty}\sup\|\phi_n(ab)-\phi_n(a)\phi_n(b)\|_{2,\mathrm{tr}}=\|\phi(ab)-\phi(a)\phi(b)\|_{2,\mathrm{tr}}=0$$
   and
   $$\lim_{n\to\infty}|\tau(a)-\mathrm{tr}_n(\phi_n(a))|=\lim_{n\to\infty}|\tau(a)-\mathrm{tr}_n(\phi(a))|=0.$$

   This means that  $(\phi_n:A\longrightarrow \mathbb{M}_n)_{n\geq1}$ is a trace preserving approximate morphism.

   (ii)$\to$(iii) Let $(\varphi_n:A\longrightarrow\mathbb{M}_n)_{n\geq1}$ be a trace-preserving approximate morphisms. Then, for any $\tau\in \mathrm{T}(A)$,
    $$\lim_{n\to\infty}\tau_n(\varphi_n(a))=\tau(a),$$
  for all $a\in A$. This implies that $\lim_n\|\varphi_n(a)\|_{2,\mathrm{tr}}<\infty$, for all $a\in A$.  Then by a same reason in the proof of \cite[Proposition 2.3]{RS2019}, for any ultrafilter $\omega$, $(\phi_n:A\longrightarrow \mathbb{M}_n)_{n\geq1}$ induce a homomorphism $\phi:A\longrightarrow\prod_\omega \mathbb{M}_n$ given by $\phi(a)=\pi_\omega\left(\left(\phi_n(a)\right)_n\right)$, for all $a\in A$.

  (iii)$\to$(i) 
  Let $\phi:A\longrightarrow\prod_\omega\mathbb{M}_n$ be a trace-preserving  homomorphism, $F\subseteq A$ be finite and $\varepsilon>0$. Assume that $\psi:A\longrightarrow\prod_{n\geq1}\mathbb{M}_n$ and $\psi_n:A\longrightarrow \mathbb{M}_n$  are similar to the proof of \cite[Proposition 2.2]{RS2019}. Then $\theta:A\longrightarrow\prod_{n\geq1}\mathbb{M}_n$ defined by $\theta(a)=\frac{1}{2}\left(\psi(a)+\psi(a^*)^*\right)$, for all $a\in A$. Let $\theta_n:A\longrightarrow\mathbb{M}_n$  be the components of $\theta$ such that $\theta(a)=\left(\theta_n(a)\right)_n$, for all $a\in A$. Then, for any $n\geq 1$, $\theta_n$ is $*$-preserving. Consider the following set
  $$\omega_F=\left\{n\in\mathbb{N}:\|\theta_n(ab)-\theta_n(a)\theta_n(b)\|_{2,\mathrm{tr}}<\varepsilon\quad\text{and}\quad|\tau(a)-
  \mathrm{tr}_n(\theta_n(a))|<\varepsilon,
  \quad\text{for all}\quad a,b\in F\right\}.$$

  As $\phi$ is a trace-preserving  homomorphism, so $\omega_F\in\omega$. This implies that $\omega$ is nonempty. Now, if let $\varphi=\theta_n$, for any fixed $n$, then we have the desire.
\end{proof}
\begin{lemma}\label{lem4}
  Let $A$ be a separable unital C$^*$-algebra. Then $\tau_A$ is an MA trace on $A$ if and only if is  MF. 
\end{lemma}
\begin{proof}
  As $\tau_A$ is an MA trace on the separable unital C$^*$-algebra $A$, by Lemma \ref{lem3}(ii), there is a trace-preserving approximate morphisms $(\varphi_n:A\longrightarrow\mathbb{M}_n)_{n\geq1}$. Consider the embedding $\mathbb{M}_n\subseteq \mathcal{Q}$ unitally in a trace preserving way. We now replace all $\mathbb{M}_n$'s with $\mathcal{Q}$ and $\mathrm{tr}_n$ by $\tau_\mathcal{Q}$. This together Lemma \ref{lem3}(iii) implies that we can get a unital trace-preserving homomorphism $\phi:A\longrightarrow\ell^\infty(\mathcal{Q})$. Then $\phi$ induces a trace-preserving homomorphism  $\phi_\omega:A\longrightarrow \mathcal{Q}_\omega$. Then by \cite[Proposition 3.3]{RS2019}, $\tau_A$ is an MF trace on $A$. The converse is clear.
\end{proof}
Let $A$ be a unital separable C$^*$-algebra, $\tau\in\mathrm{T}(A)$ and $J_\tau=\{a\in A_\omega:\|a\|_{2,\tau}=0\}$. Then  by \cite[Remark 4.7]{KR2014}, $J_\tau$ is a closed two-sided of $A_\omega$. Moreover, $J_\tau$ is a $\sigma$-ideal in $A_\omega$ \cite[Proposition 4.6]{KR2014}. Let $N$ be the weak clouser of $A$ under GNS representation of $A$ with respect to the state $\tau$ and let $\omega$ be a free ultrafilter on $\mathbb{N}$. Consider the canonical map $\Phi:A_\omega\longrightarrow N^\omega$ that is surjective \cite[Theorem 3.3]{KR2014}. Let $b\in A_\omega$ and $B$ be the separable C$^*$-subalgebra of $A_\omega$ generated by $A$ and the element $b$. As $J_\tau$ is a $\sigma$-ideal in $A_\omega$, there is a positive contraction $e\in J_\tau\cap B'$ such that $ex=x$ for all $J_A\cap B$. Then by \cite[Remark 4.7]{KR2014}, we have
\begin{equation}\label{eqbp}
  c=(1-e)b(1-e)\in A_\omega\cap A'.
\end{equation}

In the following, we generalize \cite[Proposition 3.6]{G2017}.
\begin{proposition}\label{pr1}
  Let $A$ be a separable unital C$^*$-algebra and $\tau_A$ be an MA trace on $A$. Then there is an order zero map $\varphi:A\longrightarrow \mathcal{Q}_\omega$ such that $\tau_{\mathcal{Q}_\omega}\left(\varphi(a)\varphi(1_A)^{n-1}\right)=\tau_A(a)$, for all $a\in A$.
\end{proposition}
\begin{proof}
   By Lemma \ref{lem4},  there is a  trace-preserving $\phi_\omega:A\longrightarrow \mathcal{Q}_\omega$. Similar to the proof of \cite[Proposition 3.6]{G2017}, define $J_{\mathcal{Q},\tau_\mathcal{Q}}=\{x\in\mathcal{Q}_\omega:\tau_\mathcal{Q}(x^*x)=0\}$, where $\tau_\mathcal{Q}$ is the trace that have obtained in the proof of Lemma \ref{lem4}. As $\tau_\mathcal{Q}$ is a trace, so $J_{\mathcal{Q},\tau_\mathcal{Q}}$ is a two-sided closed ideal in $\mathcal{Q}_\omega$. Similar to the above discussion $J_{\mathcal{Q},\tau_\mathcal{Q}}$ is a $\sigma$-ideal in $\mathcal{Q}_\omega$ and hence there is a positive contraction $e\in J_{\mathcal{Q},\tau_\mathcal{Q}}\cap C^*(\phi_\omega(A))'$ such that $ec=c$, for all $c\in J_{\mathcal{Q},\tau_\mathcal{Q}}\cap C^*(\phi_\omega(A))$. Similar to  \eqref{eqbp} define $\varphi:A\longrightarrow\mathcal{Q}_\omega$ as follows:
   $$\varphi(a):=(1_{\mathcal{Q}_\omega}-e)\phi_\omega(a)(1_{\mathcal{Q}_\omega}-e),$$
   for all $a\in A$. Again, by the same reason in the proof of \cite[Proposition 3.6]{G2017}, we have $\varphi$  is an order zero map and $\tau_{\mathcal{Q}_\omega}\left(\varphi(a)\varphi(1_A)^{n-1}\right)=\tau_A(a)$, for all $a\in A$.
\end{proof}
\begin{proposition}
  Let $A$ be a separable C$^*$-algebra and $J$ be a closed two-sided ideal of $A$.
  \begin{itemize}
    \item[(i)] If the exact sequence
    \begin{equation}\label{eqes}
      0\longrightarrow J\stackrel{\imath}{\longrightarrow}A\stackrel{\pi}{\longrightarrow}\frac{A}{J}\longrightarrow0
    \end{equation}
    splits and $\tau$ is an MA trace on $A$ such that $\tau|_J=0$, then $\tau$ is an MA trace on $\frac{A}{J}$.
    \item[(ii)] If $\tau$ is an MA trace on $\frac{A}{J}$, then the induced trace $\tau_A=\tau\circ\pi$ on $A$ is an MA trace on $A$.
  \end{itemize}
\end{proposition}
\begin{proof}
  (i) As \eqref{eqes} splits, so there is a c.c.p $*$-homomorphism splitting $\gamma:\frac{A}{J}\longrightarrow A$ i.e., $\pi\circ \gamma=\mathrm{id}_A$. Since $\tau$ is an MA trace on $A$, for every finite set $F\subseteq A$ and $\varepsilon>0$, there is an $n\geq1$ and $*$-linear map $\varphi:A\longrightarrow \mathbb{M}_n$ such that, for all $a\in F$,
    \begin{equation}\label{eqes1}
    \|\varphi(aa')-\varphi(a)\varphi(a')\|_{2,\mathrm{tr}}<\varepsilon,
    \end{equation}
     and
      \begin{equation}\label{eqes2}
      |\tau(a)-\mathrm{tr}_n(\varphi(a))|<\varepsilon.
      \end{equation}

      Define $\psi:\frac{A}{J}\longrightarrow \mathbb{M}_n$ by $\psi:=\varphi\circ\gamma$. Clearly, $\psi$ is a $*$-homomorphism. Then for every finite set $F'=\{a_1+J,a_2+J,\ldots,a_k+J\}\subseteq \frac{A}{J}$ and $\varepsilon>0$, there is an $n\geq1$ such that, for all $a+J,b+J\in F'$, by \eqref{eqes1} and \eqref{eqes2}, we have
      \begin{align}\label{eqes3}
      \nonumber
         \|\psi(a+J)\psi(b+J)-\psi(ab+J)\|_{2,\mathrm{tr}}& = \|\varphi\circ\gamma(a+J)\varphi\circ\gamma(b+J)-\varphi\circ\gamma(ab+J)\|_{2,\mathrm{tr}}\\
            & =\|\varphi(a)\varphi(b)-\varphi(ab)\|_{2,\mathrm{tr}}<\varepsilon,
      \end{align}
      and
      \begin{align}\label{eqes4}
        \nonumber
        |\tau(a+J)-\mathrm{tr}_n(\psi(a+J))|&=|\tau(a)-\mathrm{tr}_n(\varphi\circ\gamma(a+J))|\\
        &=|\tau(a)-\mathrm{tr}_n(\varphi(a))|<\varepsilon.
         \end{align}

         Hence, \eqref{eqes3} and \eqref{eqes4} imply that $\tau$ is an MA trace on $\frac{A}{J}$.

(ii) For any finite subset $F=\{a_1,a_2,\ldots, a_k\}$ of $A$, $F'=\{a_1+J,a_2+J,\ldots,a_k+J\}$ is a finite subset of $\frac{A}{J}$. Similar to (i), if $\varphi:\frac{A}{J}\longrightarrow\mathbb{M}_n$ be a $*$-linear map such that satisfies  \eqref{eqes3} and \eqref{eqes4}, the by defining $\Phi:A\longrightarrow\mathbb{M}_n$ by $\Phi:=\varphi\circ\pi$, we deduce that \eqref{eqes1} and \eqref{eqes2} hold.
\end{proof}
\section{Real MF algebras}
Let $A$ be a (complex) C$^*$-algebra and $\Phi$ be an involutory $*$-antiautomorphism of $A$. Then   $A_\Phi=\{a\in A: \ \Phi(a)=a^*\}$ is a real subalgebra of $A$  such that $A_\Phi\cap iA_\Phi=\{0\}$ and $A=A_\Phi+iA_\Phi$ and $A$ is called the complexification of $A_\Phi$, for more details related to this concept, we refer to \cite{pa, sta, sto}.

Let $A$ be a C$^*$-algebra and $\Phi$ be an involutory $*$-antiautomorphism of $A$. Then for every $c=a+ib\in A$, we consider $\|c\|=\|a\|+\|b\|$ and for considering  the qusidiagonality of $A_\Phi$ we replace $M_k(\mathbb{C})$ by $M_k(\mathbb{R})$. We work with the $\|\cdot\|_1$ on  $M_k(\mathbb{R})$.

\begin{theorem}\label{mtqd}
Let $A$ be a separable C$^*$-algebra and $\Phi$ be an involutory $*$-antiautomorphism of $A$. Then $A$ is MF (MA) if and only if $A_\Phi$ is MF (MA).
\end{theorem}
\begin{proof}
Let $A_\Phi$ be an MF-algebra and let $F=\{c_1=a_1+ib_1, c_2=a_2+ib_2,\ldots, c_n=a_n+ib_n\}$ be a finite subset in $A$. Set $F_r=\{a_1, a_2,\ldots, a_n,b_1, b_2,\ldots,b_n\}$ is a finite subset in $A_\Phi$. For every $\varepsilon>0$, there exist $k\geq 1$ and a $*$-real-linear  map $\varphi:A_\Phi\longrightarrow M_k(\mathbb{R})$ such that
\begin{equation}\label{eq1qd}
    \|\varphi(ab)-\varphi(a)\varphi(b)\|<\frac{\varepsilon}{4},
\quad
\text{and}\quad
   \|\varphi(a)\|>\|a\|-\frac{\varepsilon}{4},
\end{equation}
for  all $a,b\in F_r$. Now, suppose that $\varphi^c:A\longrightarrow M_k(\mathbb{C})$ is the complexification of $\varphi$ that is a $*$-linear map. Then by \eqref{eq1qd}, we have
\begin{eqnarray*}
  \|\varphi^c(c_jc_l)-\varphi(c_j)^c\varphi^c(c_l)\| &=&\|\varphi^c\left([a_ja_l-b_jb_l]+i[b_ja_l+a_jb_l]\right)-\varphi^c(a_j+ib_j)\varphi^c(a_kj+ib_j)\|  \\
   &\leq& \|\varphi(a_ja_l-b_jb_l)-\varphi(a_j)\varphi(a_l)+\varphi(b_j)\varphi(b_l)\|\\
   &&\hspace{-0.3cm} +\|i[\varphi(b_ja_l+a_jb_l)- \varphi(b_j)\varphi(a_l)-\varphi(a_j)\varphi(b_l)\|  \\
   &<& \frac{\varepsilon}{4}+\frac{\varepsilon}{4}+\frac{\varepsilon}{4}+\frac{\varepsilon}{4} \\
   &=&\varepsilon,
\end{eqnarray*}
and
\begin{eqnarray*}
 |~ \|\varphi^c(c_k)\|-\|c_k\| ~|&=&|~\|\varphi^c(a_k+ib_k)\|-\|(a_k+ib_k)\|~|  \\
   &=&|~ \|\varphi(a_k)\|-\|a_k\|+\|\varphi(b_k)\|-\|b_k\|~|\\
   &<& \frac{\varepsilon}{4}+\frac{\varepsilon}{4} \\
   &<&\varepsilon,
\end{eqnarray*}
for every $c_k, c_l\in \mathcal{F}_r$. The above inequalities imply that $A$ is MF. Define $\sigma:\mathbb{C}\longrightarrow M_2(\mathbb{R})$ by
\begin{equation}\label{eqt1}
    \sigma(a+ib)=\left(
                   \begin{array}{cc}
                     a & b \\
                    -b & a \\
                   \end{array}
                 \right)
\end{equation}
and $\rho:M_2(\mathbb{R})\longrightarrow\mathbb{C}$ by
\begin{equation}\label{eqt2}
    \rho\left(
                   \begin{array}{cc}
                     a & b \\
                    c & d \\
                   \end{array}
                 \right)=\frac{1}{2}(a+d)+\frac{1}{2}i(b-c).
\end{equation}

Then $\sigma$ and $\rho$ are c.p. maps such that $\rho\circ\sigma$ equal to the identity map \cite[Lemma 4]{ho}. It is easy to verify that $\sigma$ is a homomorphism. Moreover, for some $k\in\mathbb{N}$, $\sigma^{(k)}:M_k(\mathbb{C})\longrightarrow M_{2k}(\mathbb{R})$ becomes a homomorphism. 
  Assume that $A$ is an MF algebra. Thus, for every finite subset ${F}$ of $A_\Phi$, $\varepsilon>0$, there exists  $k\geq1$ and linear map $\varphi:A\longrightarrow M_k(\mathbb{C})$ such that the conditions (i) and (ii) in the definition of MF algebras hold for every $a\in {F}$. Define $\varphi':A\longrightarrow M_{2k}(\mathbb{R})$ by $\varphi'=\sigma^{(k)}\circ\varphi\circ\Phi\circ*$. Clearly, $\Phi\circ*$ is a real-linear homomorphism and $\theta^{(k)}$ is real-linear, so $\varphi'$ is a real-linear map. Then
\begin{eqnarray*}
  \|\varphi'(ab)-\varphi'(a)\varphi'(b)\| &=&\|\sigma^{(k)}\circ\varphi\circ\Phi\circ*(ab)-\sigma^{(k)}\circ\varphi\circ\Phi\circ*(a)\ \sigma^{(k)}\circ\varphi\circ\Phi\circ*(b)\|  \\
   &=& \|\sigma^{(k)}\circ\varphi(ab)-\sigma^{(k)}\circ\varphi(a)\ \sigma^{(k)}\circ\varphi(b)\|\\
   &=& \|\sigma^{(k)}\left(\varphi(ab)-\varphi(a)\varphi(b)\right)\| \\
   &\leq&\|\sigma^{(k)}\|\ \|\varphi(ab)-\varphi(a)\varphi(b)\|\\
   &<&\varepsilon,
\end{eqnarray*}
and
\begin{eqnarray*}
  |~\|\varphi'(a)\|-\|a\|~| &=&|~\|\theta^{(k)}\circ\varphi\circ\Phi\circ*(a)\|-\|a\|~| =|~\|\sigma^{(k)}\circ\varphi(a)\|-\|a\|~| \\
   &\leq& |~\|\sigma^{(k)}\|~\|\varphi(a)\|-\|a\|~|\leq|~\|\varphi(a)\|-\|a\|~|\\
   &<& \varepsilon,
\end{eqnarray*}
for every $a,b\in {F}$. These show that $A$ is MF. Similarly, we can prove the case MA.
\end{proof}
\section*{Data availability statement }
All data generated or analysed during this study are included in this published article (and its supplementary information files).

\end{document}